\documentclass[11pt]{article}
\usepackage[T1]{fontenc}
\usepackage{lmodern}
\usepackage[a4paper,margin=28mm]{geometry}
\usepackage{amsmath,amssymb,amsthm,booktabs}
\usepackage{microtype}
\usepackage{authblk}
\usepackage[hidelinks]{hyperref}
\newtheorem{theorem}{Theorem}
\newtheorem{lemma}[theorem]{Lemma}

\newtheorem{corollary}[theorem]{Corollary}
\theoremstyle{definition}
\newtheorem{question}[theorem]{Question}
\theoremstyle{remark}

\newcommand{\ch}{\operatorname{ch}}

\begin{document}
\title{Counterexamples to Relaxed Multiplication in List Coloring}
\author[1]{Xiaolan Hu}
\author[2]{Jiaao Li}

\affil[1]{\small School of Mathematics and Statistics and Hubei Key Laboratory of Mathematical Sciences, Central China Normal University, Wuhan 430079, China.}
\affil[2]{\small School of Mathematical Sciences and LPMC, Nankai University, Tianjin 300071, China.\linebreak
Emails: xlhu@ccnu.edu.cn; lijiaao@nankai.edu.cn.}

\date{}
\maketitle
\begin{abstract}
Jiang and Zhu [{\em J. Combin. Theory Ser. B} 137 (2019) 112--117]
asked whether every $k$-choosable graph is $(km+1:m)$-choosable
for every positive integer $m$. We show that, for all integers
$k\ge3$, $m\ge2$, and $t\ge0$ satisfying $m>8t+\epsilon$,
where $\epsilon\in\{0,1\}$ and $m\equiv\epsilon\pmod2$, there is
a graph with chromatic number and choice number both equal to $k$ that is not
$(km+t:m)$-choosable. This answers their question negatively for
every $k\ge3$ and $m\ge10$. In particular, our construction
yields a $3$-choosable graph on $91$ vertices that is not
$(6:2)$-choosable, answering a question of Dvo\v{r}\'ak,
Hu and Sereni [{\em Adv. Comb.} (2019), Paper No.~5]. More generally, for  every fixed $b\ge1$, $m\ge2$, and $t\ge0$, we construct
an $(a:b)$-choosable graph that is not $(am+t:bm)$-choosable
for all sufficiently large $a$, with an explicit bound on $a$.
Let $F(a,b)$ be the supremum of the strong fractional choice
numbers of $(a:b)$-choosable graphs. We prove that for every fixed
$b\ge1$, $F(a,b)-\frac ab$ tends to infinity as $a\to\infty$.
\end{abstract}

\smallskip
\noindent\textbf{Keywords:} list coloring; multiple list coloring;
strong fractional choice number

\section{Introduction}

All graphs in this paper are finite and simple. A list assignment $L$ on a graph $G$
assigns a finite set $L(v)$ of colors to each vertex $v$. For a positive integer $b$, an \emph{$(L:b)$-coloring}, also called a \emph{$b$-fold $L$-coloring}, 
assigns a $b$-element subset $\phi(v)$ of $L(v)$ to each
vertex $v$ such that $\phi(u)\cap\phi(v)=\varnothing$
whenever $uv\in E(G)$. For positive integers $a\ge b$, a graph $G$ is
\emph{$(a:b)$-choosable} if it has an $(L:b)$-coloring for every list
assignment $L$ with $|L(v)|=a$ at each vertex. A graph is \emph{$k$-choosable} if it is $(k:1)$-choosable. The \emph{choice number} $\ch(G)$ is the least integer $k$
such that $G$ is $k$-choosable, and  $\chi(G)$ denotes the
chromatic
number.

Erd\H{o}s, Rubin and Taylor~\cite{ERT} asked whether every
$(a:b)$-choosable graph is $(am:bm)$-choosable for every positive
integer $m$. The analogous statement for ordinary multiple coloring follows by replacing each color with $m$ distinct copies, but this argument does not apply to arbitrary list assignments.  Tuza and Voigt~\cite{TV} proved that every $2$-choosable graph is $(2m:m)$-choosable for every positive integer $m$. Dvo\v{r}\'ak, Hu and Sereni~\cite{DHS} answered the general question negatively by constructing, for every
$k\ge4$,  a $k$-choosable graph that is not $(2k:2)$-choosable.

A natural weakening of the multiplication question allows one
additional color in each list. This  question was discussed by Jiang and Zhu~\cite{JZ}
and subsequently by Li and Zhu~\cite{LZ}.

\begin{question} \cite{JZ,LZ}
\label{q1}
Is every $k$-choosable graph $(km+1:m)$-choosable for every positive integer $m$?
\end{question}

The following theorem answers Question~\ref{q1} negatively
for every $k\ge3$ and $m\ge10$.

\begin{theorem}\label{thm}
Let $k\ge3$, $m\ge2$, and $t\ge0$ be integers, and let $\epsilon\in\{0,1\}$ satisfy $m\equiv\epsilon\pmod2$. If $m>8t+\epsilon$, then there is a graph $G$ with $\chi(G)=\ch(G)=k$ that is not $(km+t:m)$-choosable.
\end{theorem}

Taking $t=0$ in Theorem~\ref{thm} shows that, for every
$k\ge3$ and every $m\ge2$, there is a $k$-choosable graph
that is not $(km:m)$-choosable.
In particular, our construction yields a $3$-choosable graph
on $91$ vertices that is not $(6:2)$-choosable, answering
a question of Dvo\v{r}\'ak, Hu and Sereni~\cite{DHS}.
Together with the result of Tuza and Voigt~\cite{TV},
this shows that, for each fixed $m\ge2$, every
$k$-choosable graph is $(km:m)$-choosable if and only if
$k\in\{1, 2\}$.

This settles the case $b=1$ of the multiplication question
of Erd\H{o}s, Rubin and Taylor~\cite{ERT}.
We next consider whether multiplication can also fail for
each prescribed initial multiplicity $b\ge2$.
Not every pair $(a,b)$ admits a counterexample.
Indeed, if $a=2b$ and $b$ is odd, then every
$(2b:b)$-choosable graph is $2$-choosable~\cite{GT},
and hence $(2bm:bm)$-choosable for every positive
integer $m$~\cite{TV}.

The next theorem establishes an explicit bound on $a$ for which
multiplication can fail even with a prescribed number
of additional colors.

\begin{theorem}\label{thm:main}
Let $a,b,m,t$ be integers with $b\ge1$, $m\ge2$, and
$t\ge0$, and let
$\eta=4b+\left\lceil\frac{2(t+2b)}{m-1}\right\rceil$.
If $a\ge(6m\eta)^{\eta+3}$,
then there is a graph $G$ that is $(a:b)$-choosable but not
$(am+t:bm)$-choosable.
\end{theorem}

Question~\ref{q1} has an equivalent formulation in terms of strong fractional choice numbers.
The \emph{fractional choice number} of a graph $G$ is the infimum of $a/b$ over all positive integers $a,b$ for which $G$ is $(a:b)$-choosable. Alon, Tuza and Voigt~\cite{ATV}
proved that this
parameter  equals the fractional chromatic number. The \emph{strong fractional choice number}, introduced by Zhu~\cite{Zhu}, requires choosability at every ratio above a threshold:
\[
 \ch_f^s(G)=\inf\{r\in\mathbb R: G\text{ is $(a:b)$-choosable whenever }
 a,b\in\mathbb N\text{ and }a/b\ge r\ge 1\}.
\]
This parameter has been studied for planar graphs~\cite{XZ,Zhu} and
series--parallel graphs~\cite{LZ}. Xu and Zhu~\cite{XZ} proved that it is rational for every finite graph. 
The possibility that it never exceeds the
choice number was discussed by Jiang and Zhu~\cite{JZ}
and Li and Zhu~\cite{LZ}, and was explicitly formulated
as a question by Xu and Zhu~\cite{XZ}.

\begin{question}
\cite{JZ,LZ,XZ}
\label{q2}
Is $\ch_f^s(G)\le\ch(G)$ for every graph $G$?
\end{question}

To see the equivalence, let $\ch_m(G)$ be the least integer $a$ such that $G$ is $(a:m)$-choosable. By~\cite[Lemma~2.1]{XZ},
\begin{equation}\label{eq}
\ch_f^s(G)=\sup_{m\ge1}\frac{\ch_m(G)-1}{m}.
\end{equation}

Consequently, $\ch_f^s(G)\le k$ if and only if $G$ is $(km+1:m)$-choosable for every positive integer $m$. A failure of $(km:m)$-choosability alone therefore does not imply $\ch_f^s(G)>k$. In contrast, the graph in Theorem~\ref{thm} satisfies $\ch_m(G)\ge km+t+1$, and hence
$\ch_f^s(G)\ge k+\frac{t}{m}$. Taking $t=1$ and $m=10$
therefore answers Question~\ref{q2} negatively.

For each integer $k\ge3$, define
\[
 F(k)=\sup\{\ch_f^s(G):G\text{ is }k\text{-choosable}\}.
\]
Theorem~\ref{thm} yields the following quantitative bound.

\begin{corollary}\label{cor:gap}
For every integer $k\ge3$, we have $F(k)\ge k+1/8$.
\end{corollary}

\begin{proof}
For each positive integer $t$, apply Theorem~\ref{thm}
with $m=8t+2$. The resulting graph $G$ satisfies
$\ch(G)=k$ and $\ch_f^s(G)\ge k+t/(8t+2)$.
Letting $t\to\infty$ proves the assertion.
\end{proof}

The graphs in this proof may depend on $t$, so
Corollary~\ref{cor:gap} does not assert that a single graph
attains the lower bound. It is natural to ask the
following question.

\begin{question}\label{q:gap}
What is the value of $F(k)$ for each $k\ge3$?
In particular, what is the order of growth of $F(k)-k$ as $k\to\infty$?
\end{question}

The following theorem establish explicit and asymptotic lower bounds
towards Question~\ref{q:gap}. Throughout the paper, $\log$ denotes the natural logarithm and
$\mathrm e$ is its base.

\begin{theorem}\label{thm:maink}
For every integer $k\ge3$,
\[
 F(k)\ge k+
 \left\lfloor\frac{\log k}{2\log\log k+8}\right\rfloor.
\]
Moreover, $F(k)\ge k+1$ for every $k\ge341$.
\end{theorem}

We next extend $F(k)$ to prescribed initial multiplicities. For
integers $a\ge b\ge1$, define
\[
 F(a,b)=\sup\{\ch_f^s(G):G\text{ is }(a:b)\text{-choosable}\}.
\]
In particular, $F(k,1)=F(k)$. 
For fixed $b$, we obtain the
following asymptotic lower bound.

\begin{theorem}\label{cor:generalized-F}
For every fixed integer $b\ge1$,
\[
 \liminf_{a\to\infty}
 \frac{(F(a,b)-a/b)\log\log a}{\log a}
 \ge\frac{1}{2b}.
\]
\end{theorem}

In communication with Rongxing Xu, we learned that Gu and Xu~\cite{GX} independently obtained related counterexamples, and they provided some planar constructions.

The rest of the paper is organized as follows. Section~\ref{sec:basic}
proves Theorem~\ref{thm}. Section~\ref{sec:general} establishes
Theorem~\ref{thm:main} for prescribed initial multiplicities.
Sections~\ref{sec:growth} and~\ref{sec:asymptotic} prove
Theorems~\ref{thm:maink} and~\ref{cor:generalized-F}, respectively.

\section{Proof of Theorem \ref{thm}}\label{sec:basic}

We first construct a $3$-choosable graph with the required
failure of multiple choosability, and then extend the
construction to every $k\ge3$. We use the following
generalization of choosability.
Let $F$ be a graph, let $c$ be a positive integer, and let
$f:V(F)\to\mathbb N$. The graph $F$ is \emph{$(f:c)$-choosable}
if it has an $(L:c)$-coloring for every list assignment $L$
satisfying $|L(v)|=f(v)$ for every $v\in V(F)$. When $c=1$,
we simply say that $F$ is \emph{$f$-choosable}.

Our construction uses the graph $H$ with
\[
 V(H)=\{u,v,w,x,y,z\},\qquad
 E(H)=\{uv,vw,uw,xy,yz,xz,ux,vy,wz\}.
\]
Thus $H$ is the triangular prism, consisting of the triangles
$uvw$ and $xyz$, together with the matching edges $ux,vy,wz$.
Define
\[
 f(u)=f(v)=f(y)=f(z)=2,\qquad f(w)=f(x)=3.
\]
Arranging the vertices of $H$ as
$\left(\begin{smallmatrix}u&v&w\\x&y&z\end{smallmatrix}\right)$,
the function $f$ corresponds to the list-size array
$\left(\begin{smallmatrix}2&2&3\\3&2&2\end{smallmatrix}\right)$.
The following result is the $2\times3$ case
proved by Isaak~\cite{Isaak}.

\begin{lemma} {\rm (\cite{Isaak})}
\label{lem:ordinary}
The graph $H$ is $f$-choosable.
\end{lemma}

We next specify lists on $H$ whose sizes are obtained by
multiplying $f$ by $m$ and adding $t$. The following lemma
determines exactly when these lists admit an $m$-fold coloring.

\begin{lemma}\label{prop:local}
Let $s\ge1$, $t\ge0$, and $\epsilon\in\{0,1\}$ be integers, and let
$m=2s+\epsilon$. Let pairwise disjoint color sets $A,B,C,D,Z$ satisfy
\[
 |A|=|B|=|C|=s,\qquad |D|=3s+t+2\epsilon,\qquad |Z|=\epsilon.
\]
Define lists on $H$ by
\[
 \begin{aligned}
 L(u)&=A\cup D,&L(z)&=B\cup D,\\
 L(v)&=L(y)=C\cup D,&L(w)&=L(x)=A\cup B\cup C\cup D\cup Z.
 \end{aligned}
\]
Then $|L(v)|=mf(v)+t$ for every $v\in V(H)$, and $H$
has an $(L:m)$-coloring if and only if $s\le4t$.
\end{lemma}

\begin{proof}
The list sizes follow directly from $m=2s+\epsilon$ and the
definition of $f$.

Suppose first that $H$ has an $(L:m)$-coloring $\phi$.
Assign weights $\omega(u)=\omega(v)=\omega(y)=\omega(z)=2$ and
$\omega(w)=\omega(x)=1$,
and let $P=A\cup B\cup C\cup D\cup Z$.
For each $\gamma\in P$, the set
\[
 I_\gamma=\{v\in V(H):\gamma\in\phi(v)\}
\]
is independent in $H$. Its weight is at most two if
$\gamma\in A\cup B\cup Z$, at most three if $\gamma\in C$,
and at most four if $\gamma\in D$.
Indeed, the vertices whose lists contain a color of $A$ or $B$
induce a two-edge star with center of weight two and leaves
of weight one. For a color of $C$, they induce the path $wvyx$,
whose maximum independent-set weight is three.
For a color of $D$, the bound follows because an independent
set meets each of the triangles $uvw$ and $xyz$ in at most
one vertex. Finally, a color of $Z$ occurs only at $w,x$.
Counting weighted color incidences yields
\[
 10m
 =\sum_{\gamma\in P}\sum_{v\in I_\gamma}\omega(v)
 \le 2s+2s+3s+4(3s+t+2\epsilon)+2\epsilon
 =19s+4t+10\epsilon.
\]
Since $10m=20s+10\epsilon$, we obtain $s\le4t$.

Conversely, suppose that $s\le4t$. First assume $\epsilon=0$
and let $q=\lceil s/4\rceil$. Empty parts are allowed in the
following partitions.
Let $A_0,A_1$ be a partition of $A$, and let $B_0,B_1$
be a partition of $B$, with
\[
 |A_0|=|B_0|=q,\qquad |A_1|=|B_1|=s-q.
\]
Let $C_1,C_2$ be a partition of $C$ with
$|C_1|=\lceil s/2\rceil$ and $|C_2|=\lfloor s/2\rfloor$.
Choose pairwise disjoint subsets
$D_{uy},D_{vz},D_{uz},D_{vx},D_{yw}$ of $D$ with
\[
 \begin{aligned}
 |D_{uy}|=|D_{vz}|&=2s-2q,&
 |D_{uz}|&=q,\\
 |D_{vx}|&=2q-\lceil s/2\rceil,&
 |D_{yw}|&=2q-\lfloor s/2\rfloor.
 \end{aligned}
\]
These sizes are nonnegative since $q\le s$ and
$\lceil s/2\rceil\le2q$. Their sum is
$3s+q\le3s+t=|D|$.
Define
\[
 \begin{aligned}
 \phi(u)&=A_0\cup D_{uy}\cup D_{uz},&
 \phi(z)&=B_0\cup D_{vz}\cup D_{uz},\\
 \phi(v)&=C_1\cup D_{vz}\cup D_{vx},&
 \phi(y)&=C_2\cup D_{uy}\cup D_{yw},\\
 \phi(w)&=A_1\cup B_1\cup C_2\cup D_{yw},&
 \phi(x)&=A_1\cup B_1\cup C_1\cup D_{vx}.
 \end{aligned}
\]
For every $v\in V(H)$, we have $\phi(v)\subseteq L(v)$
and $|\phi(v)|=2s$. Moreover, each color is assigned to
an independent set, so $\phi(u)\cap\phi(v)=\varnothing$
whenever $uv\in E(H)$. Thus $\phi$ is an $(L:2s)$-coloring.

Now assume $\epsilon=1$. Choose distinct colors
$\gamma_1,\gamma_2\in D$ and let $Z=\{\gamma_3\}$.
Remove $\gamma_1,\gamma_2,\gamma_3$ from every list in which
they occur. The resulting lists have the form already treated
for $\epsilon=0$, with $D$ replaced by
$D\setminus\{\gamma_1,\gamma_2\}$, so they admit a $2s$-fold
coloring. Assign $\gamma_1$ additionally to $u,y$,
$\gamma_2$ to $v,z$, and $\gamma_3$ to $w,x$.
These pairs are independent, and each vertex receives one
additional color from its list. Hence $H$ has an
$(L:2s+1)$-coloring.
\end{proof}

We use the following extension of the construction in the
concluding remarks of Dvo\v{r}\'ak, Hu and Sereni~\cite{DHS}
both for $k=3$ and for larger choice numbers.

\begin{lemma}\label{lem:lift}
Let $k,a,b,c$ be positive integers with $k\ge c$ and $a\ge b$.
Let $F$ be a nonempty graph and let
$f:V(F)\to\{c,\ldots,k\}$ satisfy
\[
 j(v)=\frac{k-f(v)}{c}\in\mathbb Z_{\ge0}
 \qquad\text{for every }v\in V(F).
\]
Suppose that $F$ is $(f:c)$-choosable and has no
$(L:b)$-coloring for a list assignment $L$ satisfying
\[
 |L(v)|=a-bj(v)\ge b
 \qquad\text{for every }v\in V(F).
\]
Let $d=\max_{v\in V(F)}j(v)$. Then there is a graph $G$ that is
$(k:c)$-choosable but not $(a:b)$-choosable, with
\[
 |V(G)|=d+\binom{a}{b}^{d}|V(F)|.
\]
\end{lemma}

\begin{proof}
If $d=0$, then $j(v)=0$ for every $v\in V(F)$. Hence
$f(v)=k$ and $|L(v)|=a$ for every $v\in V(F)$. Thus $F$ is
$(k:c)$-choosable but not $(a:b)$-choosable. Since
\[
 d+\binom ab^d|V(F)|=|V(F)|,
\]
the assertion holds with $G=F$. We may therefore assume $d\ge1$.

Choose pairwise disjoint $a$-element color sets
$Q_1,\ldots,Q_d$, each disjoint from
$\bigcup_{v\in V(F)}L(v)$.
For each tuple $\mathbf S=(S_1,\ldots,S_d)$ with
$S_i\in\binom{Q_i}{b}$, take a vertex-disjoint copy
$F_{\mathbf S}$ of $F$, with vertex $v_{\mathbf S}$
corresponding to $v\in V(F)$.

Add independent vertices $u_1,\ldots,u_d$, join each
$v_{\mathbf S}$ to $u_1,\ldots,u_{j(v)}$, and add no other
edges. Denote the resulting graph by $G$. Its number of
vertices is
\[
 |V(G)|=d+\binom{a}{b}^{d}|V(F)|.
\]

We first prove that $G$ is $(k:c)$-choosable. Consider any
list assignment with $k$ colors at each vertex of $G$.
Choose $c$ colors for each of $u_1,\ldots,u_d$ and remove
these colors from the lists of their neighbors. Each vertex
$v_{\mathbf S}$ retains at least
$k-cj(v)=f(v)$
colors. Since $F$ is $(f:c)$-choosable, the coloring extends
to every copy $F_{\mathbf S}$. Thus $G$ is
$(k:c)$-choosable.

Define a list assignment $L'$ on $G$ by
\[
 L'(u_i)=Q_i \quad (1\le i\le d),\qquad
 L'(v_{\mathbf S})
 =L(v)\cup\bigcup_{i=1}^{j(v)}S_i.
\]
The sets in the union are pairwise disjoint, so every list
has size $|L'(v_{\mathbf S})|
 =|L(v)|+bj(v)=a$.

Suppose that $G$ has an $(L':b)$-coloring $\phi$, and let
$\mathbf S_0=(\phi(u_1),\ldots,\phi(u_d))$.
Since $\phi(u_i)\in\binom{Q_i}{b}$ for every $i$,
the graph $F_{\mathbf S_0}$ is one of the copies in the
construction. Each vertex $v_{\mathbf S_0}$ is adjacent to
$u_1,\ldots,u_{j(v)}$, and hence cannot use any of the colors
added to $L(v)$. Therefore,
$\phi(v_{\mathbf S_0})\subseteq L(v)$ for every $v\in V(F)$.
Thus $\phi$ induces an $(L:b)$-coloring of $F$, a contradiction.
Therefore, $G$ is not $(a:b)$-choosable.
\end{proof}

When $c=1$, the disjoint union of $K_k$ and the graph obtained
in Lemma~\ref{lem:lift} is not $(a:b)$-choosable and has both
chromatic number and choice number equal to $k$. If, in addition, $f\equiv k-1$ and
$\chi(F)=k-1$, then the graph $G$ constructed in the proof
already satisfies $\chi(G)=\ch(G)=k$: in this case $d=1$,
and the new vertex is adjacent to every vertex in each copy
of $F$.

\begin{proof}[Proof of Theorem~\ref{thm}]
Let $m=2s+\epsilon$. Since $m\ge2$ and
$m>8t+\epsilon$, we have $s\ge1$ and $s>4t$.
By Lemmas~\ref{lem:ordinary} and~\ref{prop:local}, $H$ is
$f$-choosable and has a list assignment $L$ with
\[
 |L(v)|=mf(v)+t
 \qquad\text{for every }v\in V(H)
\]
that admits no $(L:m)$-coloring.

Apply Lemma~\ref{lem:lift} with
\[
 F=H,\qquad k=3,\qquad a=3m+t,\qquad b=m,\qquad c=1.
\]
Since
\[
 a-b(3-f(v))=mf(v)+t=|L(v)|,
\]
the lemma yields a $3$-choosable graph $G_3$ that is not
$(3m+t:m)$-choosable. The construction contains a copy of
$H$, which contains a triangle. Hence
$\chi(G_3)=\ch(G_3)=3$.

In this construction, $d=1$, and $G_3$ consists of one new
vertex and $\binom{3m+t}{m}$ copies of $H$. Each copy has
six vertices and nine edges, and is joined to the new vertex
by four edges. Thus
\[
 |V(G_3)|=1+6\binom{3m+t}{m},
 \qquad
 |E(G_3)|=13\binom{3m+t}{m}.
\]

Now suppose that $3\le j<k$ and that there is a graph $G_j$
with $\chi(G_j)=\ch(G_j)=j$ that is not
$(jm+t:m)$-choosable. Apply Lemma~\ref{lem:lift} with
\[
 F=G_j,\qquad f\equiv j,\qquad k=j+1,\qquad
 a=(j+1)m+t,\qquad b=m,\qquad c=1.
\]
The required list size on $G_j$ is
\[
 a-b=(j+1)m+t-m=jm+t.
\]
Hence the lemma yields a graph $G_{j+1}$ that is not
$((j+1)m+t:m)$-choosable. By the observation following
Lemma~\ref{lem:lift},
$\chi(G_{j+1})=\ch(G_{j+1})=j+1$.
Therefore Theorem~\ref{thm} follows by induction.
\end{proof}

\section{Proof of Theorem~\ref{thm:main}}\label{sec:general}

The proof of Theorem~\ref{thm:main} combines three steps. A Small Pot argument reduces ordinary
list coloring to a multicoloring bound for a Kneser graph. Replacing each
vertex by a clique then yields the required initial multiplicity.
Finally, a list assignment with a palette smaller than the total demand
provides the obstruction, to which Lemma~\ref{lem:lift} applies.

For a graph $G$ and a positive integer $b$, the $b$-fold chromatic
number $\chi_b(G)$ is the smallest number of colors needed to assign
a set of $b$ colors to each vertex of $G$ such that adjacent vertices
receive disjoint sets.
For integers $n\ge s\ge1$, the Kneser graph $KG_{n,s}$ is the graph
whose vertices are the $s$-subsets of $[n]=\{1,\ldots,n\}$, two such
subsets being adjacent if and only if their intersection is empty.
The complement of $KG_{n,s}$ is denoted by $\overline{KG}_{n,s}$.
For a positive integer $r$, let $rKG_{n,s}$ denote the disjoint union
of $r$ copies of $KG_{n,s}$.
The join $A\vee B$ of two vertex-disjoint graphs is obtained
from their disjoint union by adding all edges between $A$ and $B$.
For a graph $F$ and a positive integer $c$, the graph $F[K_c]$
is obtained by replacing each vertex $v$ of $F$ with a clique
$\{(v,1),\ldots,(v,c)\}$ and adding all edges between two such
cliques whenever their corresponding vertices are adjacent in $F$.
In particular, $|V(F[K_c])|=c|V(F)|$.

We use the following form of the Small Pot Lemma for ordinary
list coloring.

\begin{lemma}[\cite{CR}]\label{lem:smallpot}
Let $G$ be a graph and $f:V(G)\to\mathbb N$ satisfy
$f(v)<|V(G)|$ for every $v\in V(G)$.
If $G$ is not $f$-choosable, then there is a list assignment $L$
with $|L(v)|=f(v)$ for every $v\in V(G)$, using fewer than
$|V(G)|$ colors in total, such that $G$ has no $L$-coloring.
\end{lemma}

We also use a lower bound of Stahl for the multiple chromatic
numbers of Kneser graphs. 

\begin{lemma}[{\cite{Stahl}}]\label{lem:stahl}
Let $n\ge2s\ge4$, $k\ge1$, and $0\le j<s$ be integers. Then
$\chi_{sk-j}(KG_{n,s})
 \ge kn-2j-(s^2-3s+4)$.
In particular, if $d=s^2-s+3$, then
$\chi_{sk+1}(KG_{n,s})\ge nk+n-d+1$ for every integer $k\ge0$.
\end{lemma}

\begin{lemma}\label{lem:multimonotone}
For every nonempty graph $G$ and integers $t\ge1$ and $s\ge0$,
$\chi_{t+s}(G)\ge\chi_t(G)+s$.
\end{lemma}

\begin{proof}
It suffices to prove the assertion for $s=1$.
Consider a $(t+1)$-fold coloring using $\chi_{t+1}(G)$ colors
and delete one color from its palette. Each vertex retains at
least $t$ colors. Truncating  the remaining sets to size $t$ yields a
$t$-fold coloring using at most $\chi_{t+1}(G)-1$ colors.
Thus $\chi_t(G)\le\chi_{t+1}(G)-1$.
\end{proof}

\begin{lemma}\label{lem:pairing}
Let $H$ be a nonempty graph, and let $q,k,p,h$ be positive
integers satisfying
$p+q-k>0$, $p>h$, $h-q\ge k$,
and $\chi_{h-q+1}(H)\ge p+1$.
Let $F$ be a graph on $k$ vertices that contains $q$
pairwise vertex-disjoint induced subgraphs, each isomorphic
to $\overline H$. Let
$J=K_{p+q-k}\vee F$,
and define $f:V(J)\to\mathbb N$ by
\[
 f(v)=
 \begin{cases}
  p, & v\in V(K_{p+q-k}),\\
  h, & v\in V(F).
 \end{cases}
\]
Then $J$ is $f$-choosable.
\end{lemma}

\begin{proof}
Suppose to the contrary that $J$ is not $f$-choosable.
Since $|V(J)|=p+q$ and $h<p$, we have
$f(v)<|V(J)|$ for every vertex $v$. By Lemma~\ref{lem:smallpot},
there is a list assignment $L$ with $|L(v)|=f(v)$ for every
$v\in V(J)$, using at most $p+q-1$ colors in total, such that
$J$ has no $L$-coloring.

Let $F_1,\ldots,F_q$ be pairwise vertex-disjoint induced
subgraphs of $F$, each isomorphic to $\overline H$. We
successively choose two nonadjacent vertices with a common
list color in each $F_i$, using distinct colors for different
pairs.

Suppose that such pairs have been chosen in
$F_1,\ldots,F_i$, where $0\le i<q$. Delete the $2i$ chosen
vertices and remove the $i$ used colors from all remaining
lists. Every list in $F_{i+1}$ then has at least $h-i$ colors,
and the union of the remaining lists contains at most
$p+q-1-i$ colors.
By Lemma~\ref{lem:multimonotone} and the assumption on $H$,
$\chi_{h-i}(H)
 \ge \chi_{h-q+1}(H)+q-1-i
 \ge p+q-i$.
If no two nonadjacent vertices in $F_{i+1}$ have a common
list color, then the lists of adjacent vertices in $H$ are
disjoint. Choosing $h-i$ colors from each list would produce an
$(h-i)$-fold coloring of $H$ using at most $p+q-1-i$
colors, a contradiction.

Thus two nonadjacent vertices in $F_{i+1}$ have a common
list color. Assign this color to both vertices. Since the
previously used colors have been removed, the new color is
distinct from them. Repeating this argument, we obtain $q$ pairs
colored with $q$ distinct colors.

Delete the $2q$ colored vertices and remove the $q$ used
colors from the remaining lists. Each of the $k-2q$ remaining
vertices of $F$ has at least
$h-q\ge k$
available colors, and hence these vertices can be colored
greedily. At most
$q+(k-2q)=k-q$
colors are used on $F$.
Every list on $K_{p+q-k}$ therefore retains at least
$p-(k-q)=p+q-k$
colors. The clique can consequently be colored greedily.
This is an $L$-coloring of $J$, a contradiction.
\end{proof}

\begin{proof}[Proof of Theorem~\ref{thm:main}]

Let
\[
 z=\left\lceil\frac{t+bm-1}{bm}\right\rceil,\qquad
 s=\left\lfloor\frac{2bmz-1}{m-1}\right\rfloor,\qquad
 w=2bmz-(m-1)s,
\]
and define
\[
 \begin{gathered}
 d=s^2-s+3,\qquad
 c=\left\lceil\frac{d-s+bz}{b}\right\rceil,\\
 n=s+bmc,\qquad \ell=\binom ns,\qquad k=bz\ell.
 \end{gathered}
\]
Then we have
\[
 z\ge1,\qquad 1\le w\le m-1,\qquad
 s\ge2bz\ge2,\qquad n\ge2s.
\]
For the last inequality, since $ d-2s+bz=s^2-3s+3+bz>0$ for $s\ge2$, we have $bc\ge d-s+bz>s$, and hence
$n=s+bmc>2s$.

Choose $r\in\{0,\ldots,b-1\}$ with $a\equiv r\pmod b$, and
define
\[
 q=bz-r,\qquad h=sk+bz+r,\qquad
 p=nk+bz+r+b(m-1)c.
\]
Then $1\le q\le bz$, $bz+r<s$, and
\begin{equation}\label{eq:general-congruences}
 p\equiv h\equiv r\pmod b,\qquad
 p-h=bc(mk+m-1)>0.
\end{equation}

We first verify that $p\le a$. By the definitions of $z$,
$s$, and $\eta$,
\[
 s\le\frac{2t+4bm}{m-1}
 =4b+\frac{2(t+2b)}{m-1}\le\eta.
\]
Since $bz+b-1\le s-1$, the definition of $c$ implies
$bc\le d-1$. Hence
$n\le s+m(d-1)\le2ms^2$.
Moreover,
$bz+r+b(m-1)c<n$,
and $k+1=bz\ell+1\le s\ell$.

Using $\ell=\binom ns\le
 \left(\frac{\mathrm e n}{s}\right)^s$,
we obtain
$$p<n(k+1)\le ns\ell \le ns\left(\frac{\mathrm e n}{s}\right)^s\le2ms^3(2\mathrm e ms)^s\le(6ms)^{s+3}\le(6m\eta)^{\eta+3}\le a.$$

Since 
$p+q-k=(n-1)k+2bz+b(m-1)c>0$, and $p+q\equiv k\equiv0\pmod b$, the number $(p+q-k)/b$ is a positive integer.
Let
\[
 F=K_{(p+q-k)/b}\vee\overline{zKG_{n,s}}.
\]
The graph $\overline{zKG_{n,s}}$ is the join of  $z$ 
 copies of $\overline{KG}_{n,s}$. 
Define $f:V(F)\to\mathbb N$ by
\[
 f(v)=
 \begin{cases}
  p, & v\in V(K_{(p+q-k)/b}),\\
  h, & v\in V(\overline{zKG_{n,s}}).
 \end{cases}
\]

We prove that $F$ is $(f:b)$-choosable by considering
$J=F[K_b]$. Let
$W=\bigl(\overline{zKG_{n,s}}\bigr)[K_b]$.
Then
\[
 J=K_{p+q-k}\vee W,
 \qquad |V(W)|=bz\ell=k.
\]
For each vertex $v$ of $F$, define $f_J$ to be $f(v)$ on
every vertex of the replacement clique corresponding to $v$.

The graph $W$ contains $bz$ pairwise vertex-disjoint induced
subgraphs isomorphic to $\overline{KG}_{n,s}$. Indeed, the
vertices in the clique expansion of each of the $z$ copies of
$\overline{KG}_{n,s}$ can be partitioned into $b$ sets, each
containing one vertex from every replacement clique and
inducing a copy of $\overline{KG}_{n,s}$. Since $1\le q\le bz$,
we may select $q$ of these subgraphs.

Recall that $p+q-k>0$ and $h<p$. Since
$h-q=sk+2r\ge k$,
by Lemmas~\ref{lem:stahl} and~\ref{lem:multimonotone}, we have
$$\chi_{h-q+1}(KG_{n,s})
 \ge\chi_{sk+1}(KG_{n,s})+2r
 \ge nk+n-d+1+2r\ge p+1.$$
The last inequality follows from $bc\ge d-s+bz$ and
\[
 nk+n-d+2r-p=s+bc-d-bz+r\ge0.
\]
Therefore, Lemma~\ref{lem:pairing}, applied with
$H=KG_{n,s}$ and $W$ in place of its graph $F$, shows that
$J$ is $f_J$-choosable.

For any list assignment $L$ on $F$ with $|L(v)|=f(v)$ for
every $v\in V(F)$,  assign $L(v)$ to each vertex of the
replacement clique corresponding to $v$ in $J$.
An ordinary coloring of $J$ assigns $b$ distinct colors to
each replacement clique and induces an $(L:b)$-coloring of
$F$. Thus $F$ is $(f:b)$-choosable.

Let $t'=mq-w$. By the definitions of $q$ and $w$,
$t'=m(bz-r)-\bigl(2bmz-(m-1)s\bigr)
    =(m-1)s-m(bz+r)$. Since the definition of $z$ implies
$z-1\ge\frac{t-1}{bm}$, and since $r\le b-1$ and $w\le m-1$, 
we have
$$t'\ge m(bz-b+1)-(m-1)
 =bm(z-1)+1\ge t.$$
Choose pairwise disjoint color sets $A_1,\ldots,A_n$, each of
size $mk+m-1$. In each copy of $\overline{KG}_{n,s}$, assign
the list
$\bigcup_{i\in X}A_i$
to the vertex indexed by the $s$-subset $X$ of $[n]$. Assign
the full palette $\bigcup_{i=1}^n A_i$ to every vertex of
$K_{(p+q-k)/b}$.

These lists have size $mf(v)+t'$ at every vertex $v$. Indeed, by the definitions of $h$ and $t'$,
$$s(mk+m-1)=m(sk+bz+r)+((m-1)s-m(bz+r))=mh+t'.$$
 Similarly, using $n=s+bmc$ and the definition of $p$, we have  $$n(mk+m-1)=m(nk+bz+r+b(m-1)c)+((m-1)s-m(bz+r))=mp+t'.$$

For each color in $A_i$, the vertices whose lists contain
that color form a clique: within each copy their indices
contain $i$, the copies are pairwise joined, and each is
joined to $K_{(p+q-k)/b}$. Thus each color can be used at
most once in a $bm$-fold coloring. Such a coloring would
require
\[
 bm|V(F)|=m(p+q)=mp+t'+w>mp+t'
\]
distinct colors, whereas the palette has only $mp+t'$ colors.
Hence these lists admit no $bm$-fold coloring.

Since $t'\ge t$, we may truncate the lists to obtain a list
assignment $L_0$ with
\[
 |L_0(v)|=mf(v)+t\qquad\text{for every }v\in V(F),
\]
which still admits no $(L_0:bm)$-coloring.

By~\eqref{eq:general-congruences} and $a\ge p>h\ge b$,
\[
 \frac{a-f(v)}b\in\mathbb Z_{\ge0}
 \qquad\text{and}\qquad
 |L_0(v)|
 =am+t-bm\frac{a-f(v)}b\ge bm
\]
for every $v\in V(F)$. Applying Lemma~\ref{lem:lift} to $F$
and $L_0$, we obtain an $(a:b)$-choosable graph that is not
$(am+t:bm)$-choosable.
\end{proof}

\section{Proof of Theorem~\ref{thm:maink}}\label{sec:growth}

We first prove the explicit bound $F(k)\ge k+1$ for $k\ge341$.
We then present a family of constructions with larger additive gaps and
estimate the choice numbers required by these
constructions.

\begin{lemma}\label{lem:kneser12}
The graph $KG_{12,3}$ has $220$ vertices and satisfies
\[
 \chi_{85}(KG_{12,3})=344,
 \qquad
 \ch(\overline{KG}_{12,3})\le82.
\]
\end{lemma}

\begin{proof}
By Stahl's theorem~\cite{Stahl}, for integers $n\ge6$,
$q\ge1$, and $0\le r\le2$,
$\chi_{3q-r}(KG_{n,3})=qn-2r$.
For $n=12$, $q=29$, and $r=2$, we obtain
$\chi_{85}(KG_{12,3})=344$.

By Baranyai's theorem~\cite{Baranyai}, the triples of set $[12]$
can be partitioned into $55$ classes, each consisting of four
pairwise disjoint triples. Each class is independent in
$\overline{KG}_{12,3}$. Hence $\overline{KG}_{12,3}$ is a
subgraph of the complete $55$-partite graph with parts of size
four. Kierstead, Salmon and Wang~\cite{KSW} proved that the complete
$t$-partite graph with parts of size four has choice number
$\lceil(3t-1)/2\rceil$. Therefore
$\ch(\overline{KG}_{12,3})
 \le \left\lceil\frac{3\cdot55-1}{2}\right\rceil=82$.
\end{proof}
\begin{lemma}\label{lem:unit}
For every integer $k\ge341$ and every positive integer
$m\equiv2\pmod3$, there is a graph $G$ with
$\chi(G)=\ch(G)=k$ that is not $((k+1)m-1:m)$-choosable.
Consequently, $F(k)\ge k+1$ for every $k\ge341$.
\end{lemma}

\begin{proof}
By Lemma~\ref{lem:kneser12}, the graph $KG_{12,3}$ has
$220$ vertices and satisfies
\[
 \chi_{85}(KG_{12,3})=344,\qquad
 \ch(\overline{KG}_{12,3})\le82.
\]

Let $H=K_{122}\vee\overline{KG}_{12,3}$. Then
$|V(H)|=122+\binom{12}{3}=342$. Each vertex of
$\overline{KG}_{12,3}$ is a triple of $[12]$. Define
$f:V(H)\to\mathbb N$ by
\[
 f(v)=
 \begin{cases}
  341, & v\in V(K_{122}),\\
  83, & v\in V(\overline{KG}_{12,3})\text{ and }1\in v,\\
  85, & v\in V(\overline{KG}_{12,3})\text{ and }1\notin v.
 \end{cases}
\]
There are $\binom{11}{2}=55$ vertices with $f(v)=83$ and
$\binom{11}{3}=165$ vertices with $f(v)=85$.

We prove that $H$ is $f$-choosable. Suppose to the contrary
that it is not. Since $f(v)<|V(H)|$ for every vertex $v$,
Lemma~\ref{lem:smallpot} implies that there is a list assignment
$L$ with $|L(v)|=f(v)$ for every vertex $v$ and at most $341$
colors in total, such that $H$ has no $L$-coloring.

Two vertices indexed by disjoint triples have a common list color.
Otherwise, the lists of every two adjacent vertices of
$KG_{12,3}$ are disjoint. Add two new colors to the lists of all
vertices indexed by triples containing $1$. These vertices form an
independent set in $KG_{12,3}$, and every resulting list has size
$85$. Thus the lists form an $85$-fold coloring of $KG_{12,3}$
with at most $343$ colors, contradicting
$\chi_{85}(KG_{12,3})=344$.

Color the two nonadjacent vertices of $\overline{KG}_{12,3}$ with
their common list color. Delete these vertices and remove their
color from the remaining lists on $\overline{KG}_{12,3}$. Every
remaining list has size at least $82$. Since
$\ch(\overline{KG}_{12,3})\le82$, the coloring extends to
$\overline{KG}_{12,3}$.
At most $219$ colors are used on $\overline{KG}_{12,3}$.
Every list on $K_{122}$ therefore retains at least
$341-219=122$ colors, and the clique can be colored greedily.
This is an $L$-coloring of $H$, contradicting the choice of $L$.
Hence $H$ is $f$-choosable.

Let $m\equiv2\pmod3$, $t=m-1$, and $b=(86m-1)/3$.
Choose pairwise disjoint sets $A_1,\ldots,A_{12},D$ with
\[
 |A_1|=b-2m,\qquad |A_i|=b\quad(2\le i\le12),
 \qquad |D|=3.
\]
For each triple $X$, assign the list $\bigcup_{i\in X}A_i$ to
the corresponding vertex of $\overline{KG}_{12,3}$. Assign the
full palette $D\cup\bigcup_{i=1}^{12}A_i$ to every vertex of
$K_{122}$. These lists have size $mf(v)+t$: the size is
$3b-2m=83m+t$ when $1\in X$, and $3b=85m+t$ when $1\notin X$.
The list at each clique vertex has size
\[
 |D|+\sum_{i=1}^{12}|A_i|=3+12b-2m=342m-1=341m+t.
\]

For every color, the vertices whose lists contain that color form
a clique in $H$: the triples
containing a fixed element form a clique in the complement, and
the clique $K_{122}$ is joined to it. Thus an $m$-fold
coloring would require $342m$ distinct colors, but the palette
has only $342m-1$ colors. The lists admit no $m$-fold coloring.

For every $k\ge341$, Lemma~\ref{lem:lift} applies with
$a=km+t$, $b=m$ and $c=1$. Taking the disjoint union of the
resulting graph with $K_k$, we obtain a graph $G$ with
$\chi(G)=\ch(G)=k$ that is not
$((k+1)m-1:m)$-choosable. As $m$ tends to infinity through
integers congruent to $2$ modulo $3$, equation~\eqref{eq} yields
$F(k)\ge k+1$.
\end{proof}

The next lemma uses Stahl's general bound.

\begin{lemma}\label{lem:growth}
Let $r\ge2$, $s>2r$, and $n\ge2s$ be integers. Define
\[
 d=s^2-s+3,\quad \ell=\binom ns,\quad
 h=s\ell+r,\quad p=n(\ell+1)-d.
\]
If $n(s-2r)>s(d-r)$, then $F(k)\ge k+r$ for every $k\ge p$.
\end{lemma}

\begin{proof}
By Lemma~\ref{lem:stahl}, we have
\[
 \chi_{h-r+1}(KG_{n,s})
 =\chi_{s\ell+1}(KG_{n,s})
 \ge n(\ell+1)-d+1=p+1.
\]
Since $n\ge2s$, we have $\ell\ge s$. Moreover,
$d+r\le s^2$, because $s\ge2r+1\ge r+3$. Therefore,
\[
 p-h=(n-s)\ell+n-d-r\ge n>0
\]
and
\[
 p+r-r\ell=(p-h)+(s-r)\ell+2r>0.
\]

Let
\[
 G=K_{p+r-r\ell}\vee\overline{rKG_{n,s}}.
\]
The graph $\overline{rKG_{n,s}}$ consists of $r$ pairwise
joined copies of $\overline{KG}_{n,s}$. Define
$f:V(G)\to\mathbb N$ by
\[
 f(v)=
 \begin{cases}
  p, & v\in V(K_{p+r-r\ell}),\\
  h, & v\in V(\overline{rKG_{n,s}}).
 \end{cases}
\]

The graph $\overline{rKG_{n,s}}$ has $r\ell$ vertices and
contains $r$ pairwise vertex-disjoint induced subgraphs
isomorphic to $\overline{KG}_{n,s}$. Moreover, $h-r=s\ell\ge r\ell$.
Together with
$p+r-r\ell>0$,  $h<p$, 
and  $\chi_{h-r+1}(KG_{n,s})\ge p+1$,
Lemma~\ref{lem:pairing} shows that $G$ is $f$-choosable.

Let $m$ be a positive multiple of $s$, and define
\[
 t=rm-s,\qquad b=\frac{hm+t}{s}.
\]
Since $h=s\ell+r$, we have
\[
 b=\ell m+\frac{2rm}{s}-1.
\]
Thus $b$ is a positive integer. Moreover, $t\ge0$, and
\[
 pm+t-nb
 =n-s+\frac{m}{s}\bigl[n(s-2r)-s(d-r)\bigr]>0
\]
by the assumption of the lemma.

Choose pairwise disjoint color sets $A_1,\ldots,A_n,B$ with
\[
 |A_i|=b\quad(1\le i\le n),
 \qquad
 |B|=pm+t-nb.
\]
For each $s$-subset $X$ of $[n]$, assign the list
$\bigcup_{i\in X}A_i$ to the vertex indexed by $X$ in each
copy of $\overline{KG}_{n,s}$. Assign the full palette
$B\cup\bigcup_{i=1}^n A_i$ to every vertex of
$K_{p+r-r\ell}$.
These lists have size $sb=hm+t$ on
$\overline{rKG_{n,s}}$ and size $pm+t$ on the clique.
Hence $|L(v)|=mf(v)+t$ for every vertex $v$.

For every color, the vertices whose lists contain that color
form a clique in $G$. Indeed, in each copy of
$\overline{KG}_{n,s}$, the vertices indexed by sets containing
a fixed element form a clique; the copies are pairwise joined
and are joined to $K_{p+r-r\ell}$. A color in $B$ occurs only
on this clique.
Consequently, each color can be used at most once in an
$(L:m)$-coloring. Such a coloring would require $(p+r)m$
distinct colors, whereas the palette contains only
\[
 pm+t=(p+r)m-s<(p+r)m
\]
colors. Hence $G$ has no $(L:m)$-coloring.

For every $k\ge p$, we have $f(v)\le p\le k$ and
\[
 |L(v)|=mf(v)+t=(km+t)-m(k-f(v))\ge m.
\]
Applying Lemma~\ref{lem:lift} with $c=1$, $a=km+t$, and
$b=m$ yields a $k$-choosable graph that is not
$(km+t:m)$-choosable. Therefore,
\[
 F(k)\ge k+\frac{t}{m}
 =k+r-\frac{s}{m}.
\]
Letting $m\to\infty$ through multiples of $s$ proves that
$F(k)\ge k+r$.
\end{proof}

\begin{proof}[Proof of Theorem~\ref{thm:maink}]
The bound $F(k)\ge k+1$ for every $k\ge341$ follows from
Lemma~\ref{lem:unit}. It remains to prove the first bound.
Let $\lambda=\log k$ and
\[
 r=\left\lfloor\frac{\lambda}{2\log\lambda+8}\right\rfloor.
\]
Since $k\ge3$, we have $\lambda>1$. If $r=0$, then
$F(k)\ge\ch_f^s(K_k)=k$, as required.
If $r=1$, then the definition of $r$ implies that
\[
 \lambda\ge2\log\lambda+8>8.
\]
Thus $k>\mathrm e^8>341$, and hence $F(k)\ge k+1=k+r$ by Lemma~\ref{lem:unit}.

Now suppose that $r\ge2$. Let $\rho=\log r$ and
\[
 c=\left\lceil\frac{2r}{\rho+2}\right\rceil,\qquad
 s=2r+c,\qquad d=s^2-s+3,\qquad
 n=\left\lfloor\frac{s(d-r)}{c}\right\rfloor+1.
\]
Since $s\ge2r+1\ge5$, we have
$d-r-2c=s^2-3s+3r+3>0$. It follows that $n>2s$. Moreover,
$nc>s(d-r)$. Let $\ell=\binom ns$ and $p=n(\ell+1)-d$. Since $\ell\ge s$
and $n>2s$, we have $p>0$. By Lemma~\ref{lem:growth}, we have $F(k)\ge k+r$ whenever
$k\ge p$.
It therefore remains to show that $p\le k$.

By the definition of $r$,
\[
 \lambda\ge r(2\log\lambda+8).
\]
In particular, $\lambda\ge r$, and hence
\[
 \lambda\ge r(2\rho+8).
\]
Taking logarithms, we obtain
\[
 \log\lambda\ge\rho+\log(2\rho+8).
\]
Therefore,
\[
 \log k=\lambda
 \ge r(2\log\lambda+8)
 \ge r\bigl(2\rho+2\log(2\rho+8)+8\bigr).
\]
Thus it suffices to prove that
\begin{equation}\label{eq:pbound}
 \frac{\log p}{r}
 <2\rho+2\log(2\rho+8)+8.
\end{equation}

For $r=2$, direct substitution shows that $s=6$, $n=94$, and
$p<10^{11}$. For $r=3$, we obtain $s=8$, $n=225$, and
$p<10^{17}$. For $r=4$, we obtain $s=11$, $n=400$, and
$p<10^{24}$. Hence $p<(200000r^2)^r$ for $r=2,3,4$.
Since $\rho\ge\log2$ and
$200000<\mathrm e^8(2\log2+8)^2$,
inequality~\eqref{eq:pbound} follows in these cases.

Suppose that $r\ge5$. By the definition of $c$,
\[
 \frac{2r}{\rho+2}\le c
 \le\frac{2r}{\rho+2}+1.
\]
In particular, $c\le r$ and $s\le3r$. Moreover, $\frac{4r^2}{c}\leq 2r(\rho+2)$. 
Since $s=2r+c$, we also have 
$s/r\le2+2/(\rho+2)+1/r$. Since $p=n(\ell+1)-d<2n\ell$ and
$\ell=\binom ns\le(\mathrm e n/s)^s$, we have
\[
 \frac{\log p}{r}
 <\frac{s}{r}\log\left(\frac{\mathrm e n}{s}\right)
  +\frac{\log(2n)}{r}.
\]
  We next estimate the two logarithms.

Since $s>c$ and
$s+r-3\ge1$, the definitions of $n$ and $d$ imply that
\[
 \frac ns
 \le\frac{d-r}{c}+\frac1s
 =\frac{s^2}{c}-\frac{s+r-3}{c}+\frac1s
 <\frac{s^2}{c}.
\]
Consequently,
\[
 \begin{aligned}
 \frac ns
 &<\frac{s^2}{c}
 =\frac{4r^2}{c}+4r+c\\
 &\le2r(\rho+4)+\frac{2r}{\rho+2}+1
 <2r\left(\rho+\frac92\right)=r(2\rho+9).
 \end{aligned}
\]
The last inequality follows from $\rho>1$ and $r\ge5$.
Thus
\[
 \log\left(\frac{\mathrm e n}{s}\right)
 <1+\rho+\log(2\rho+9).
\]

For $\rho>1$, we have
\[
 1+\rho+\log(2\rho+9)
 \le\frac32(\rho+2).
\]
Indeed, the difference between the right and left sides is
increasing for $\rho\ge1$ and is positive at $\rho=1$, since
$\log11<5/2$.

Using $s\le3r$ and $\rho\le r/2$, we obtain
\[
 2n<12r^2\left(\rho+\frac92\right)
 \le18r^3<\mathrm e^3r^3.
\]
Therefore $\log(2n)<3\rho+3$.

We also need the estimate
$\frac{\frac92\rho+6}{r}<\frac{14}{5}$.
To see this, note that $(\frac92\log r+6)/r$ is decreasing
for $r\ge1$. Since $r\ge5$ and $\log5<7/4$, we have
\[
 \frac{\frac92\rho+6}{r}
 \le\frac{\frac92\log5+6}{5}
 <\frac{111}{40}<\frac{14}{5}.
\]

Combining
the preceding estimates, we obtain
\[
 \begin{aligned}
 \frac{\log p}{r}
 &<\frac sr\log\left(\frac{\mathrm e n}{s}\right)
   +\frac{\log(2n)}r\\
 &\le
   \left(2+\frac{2}{\rho+2}+\frac1r\right)
   \bigl(1+\rho+\log(2\rho+9)\bigr)
   +\frac{3\rho+3}{r}\\
 &=2\bigl(1+\rho+\log(2\rho+9)\bigr)
   +\left(\frac{2}{\rho+2}+\frac1r\right)
   \bigl(1+\rho+\log(2\rho+9)\bigr)
   +\frac{3\rho+3}{r}\\
 &\le
   2\bigl(1+\rho+\log(2\rho+9)\bigr)
   +\frac{2}{\rho+2}\cdot\frac32(\rho+2)
   +\frac1r\cdot\frac32(\rho+2)
   +\frac{3\rho+3}{r}\\
 &=2\bigl(1+\rho+\log(2\rho+9)\bigr)
   +3+\frac{\frac92\rho+6}{r}\\
 &<2\bigl(1+\rho+\log(2\rho+9)\bigr)
   +3+\frac{14}{5}\\
 &=2\rho+2\log(2\rho+9)+\frac{39}{5}.
 \end{aligned}
\]
Finally, since $\rho>1$ and $\log(1+x)<x$ for $x>0$,
$$ 2\log(2\rho+9)
 =2\log(2\rho+8)
   +2\log\left(1+\frac1{2\rho+8}\right)<2\log(2\rho+8)+\frac15.$$
This proves~\eqref{eq:pbound} for $r\ge5$, and hence for every
$r\ge2$.

It follows that $\log p<\log k$, so $p<k$.
Lemma~\ref{lem:growth} therefore implies that $F(k)\ge k+r$.
Together with the cases $r=0,1$, this proves that
\[
 F(k)\ge k+
 \left\lfloor\frac{\log k}{2\log\log k+8}\right\rfloor.
\]
\end{proof}

\section{Proof of Theorem~\ref{cor:generalized-F}}\label{sec:asymptotic}

\begin{proof}
Fix an integer $b\ge1$ and a real number $\epsilon$ with
$0<\epsilon<1$. For each sufficiently large integer $a$, let
\[
 m=\left\lceil\log\log a\right\rceil,\qquad
 \eta=\left\lfloor
 \frac{(1-\epsilon)\log a}{\log\log a}
 \right\rfloor,
\]
and let
\[
 t=\left\lfloor
 \frac{(m-1)(\eta-4b)}2
 \right\rfloor-2b.
\]
Since $m,\eta\to\infty$ as $a\to\infty$, we have $m\ge3$
and $t\ge0$ for all sufficiently large $a$.

We first verify that the chosen value of $\eta$ agrees with
the definition in Theorem~\ref{thm:main}. By the definition
of $t$,
\[
 \eta-4b-\frac{2}{m-1}
 <
 \frac{2(t+2b)}{m-1}
 \le \eta-4b.
\]
Since $m\ge3$ and $\eta-4b$ is an integer, it follows that
\[
 \left\lceil\frac{2(t+2b)}{m-1}\right\rceil
 =\eta-4b.
\]
Thus $\eta
 =4b+\left\lceil\frac{2(t+2b)}{m-1}\right\rceil$,
as required.

Next, as $a\to\infty$,
\[
 m\sim\log\log a,
 \qquad
 \eta\sim\frac{(1-\epsilon)\log a}{\log\log a}.
\]
Consequently,
\[
 m\eta\sim(1-\epsilon)\log a,
 \qquad
 \log(6m\eta)=\log\log a+O(1),
\]
and hence
\[
 \frac{(\eta+3)\log(6m\eta)}{\log a}
 \longrightarrow1-\epsilon<1.
\]
Therefore, for all sufficiently large $a$,
\[
 a\ge(6m\eta)^{\eta+3}.
\]

By Theorem~\ref{thm:main}, there is an $(a:b)$-choosable
graph $G$ that is not $(am+t:bm)$-choosable. Hence
$\ch_{bm}(G)\ge am+t+1$. By~\eqref{eq},
\[
 F(a,b)\ge\ch_f^s(G)
 \ge\frac{\ch_{bm}(G)-1}{bm}
 \ge\frac ab+\frac{t}{bm}.
\]
By the definition of $t$,
\[
 t=\frac{(m-1)(\eta-4b)}2-2b+O(1).
\]
It follows that
\[
 \frac{t\log\log a}{bm\log a}
 \longrightarrow\frac{1-\epsilon}{2b}.
\]
Consequently,
\[
 \liminf_{a\to\infty}
 \frac{(F(a,b)-a/b)\log\log a}{\log a}
 \ge\frac{1-\epsilon}{2b}.
\]
Since $\epsilon>0$ is arbitrary, the result follows.
\end{proof}

\section*{Declaration of AI usage}

During the preparation of this manuscript, the authors used ChatGPT 6  to explore potential approaches and refine the presentation. All AI-assisted suggestions and arguments were carefully reviewed, revised, and independently verified by the authors. The authors take full responsibility for the mathematical correctness and content of the manuscript.

\end{document}